\documentclass[11pt]{article}
\usepackage{amsmath, amssymb, amsthm, amsbsy, ifthen, graphicx}
\date{}
\title{Variations on Cops and Robbers}
\author{
Alan Frieze \thanks{Department of Mathematical
Sciences, Carnegie Mellon University, Pittsburgh, PA 15213, email:
alan@random.math.cmu.edu. Research supported in part by NSF award
DMS-0753472.}
\and
Michael Krivelevich \thanks{School of
Mathematical Sciences, Raymond and Beverly Sackler Faculty of Exact
Sciences, Tel Aviv University, Tel Aviv 69978, Israel, e-mail:
krivelev@post.tau.ac.il. Research supported in part by USA-Israel
BSF grant 2006322, by grant 1063/08 from the Israel Science
Foundation, and by a Pazy memorial award.}
\and
Po-Shen Loh
\thanks{Department of Mathematical Sciences, Carnegie Mellon University,
Pittsburgh, PA 15217, e-mail: ploh@cmu.edu.}
}

\newtheorem{theorem}{Theorem}[section]
\newtheorem{definition}[theorem]{Definition}
\newtheorem{fact}[theorem]{Fact}
\newtheorem{inequality}[theorem]{Inequality}
\newtheorem{proposition}[theorem]{Proposition}
\newtheorem{lemma}[theorem]{Lemma}
\newtheorem{corollary}[theorem]{Corollary}

\newcommand{\pr}[1]{\mathbb{P}\left[#1\right]}

\newcommand{\bin}[1]{\text{\rm Bin}\left[#1\right]}
\newcommand{\whp}{\textbf{whp}}
\newcommand{\Gnp}{G_{n,p}}
\newcommand{\andsp}{\quad\quad\quad}
\newcommand{\del}{\partial}

\begin{document}
\maketitle

\begin{abstract}
  We consider several variants of the classical Cops and Robbers game.
  We treat the version where the robber can move $R > 1$ edges at a
  time, establishing a general upper bound of
  $n/\alpha^{(1-o(1))\sqrt{\log_\alpha n}}$, where
  $\alpha=1+\frac{1}{R}$, thus generalizing the best known upper bound
  for the classical case $R=1$ due to Lu and Peng.  We also show that
  in this case, the cop number of an $n$-vertex graph can be as large
  as $n^{1-\frac{1}{R-2}}$ for finite $R$, but linear in $n$ if $R$ is
  infinite.  For $R=1$, we study the directed graph version of the
  problem, and show that the cop number of any strongly connected
  digraph on $n$ vertices is at most $O(n(\log\log n)^2/\log n)$.  Our
  approach is based on expansion.
\end{abstract}

\section{Introduction}

The game of {\em Cops and Robbers}, introduced by Nowakowski and
Winkler \cite{NW} and independently by Quillot \cite{Qu}, is a perfect
information game played on a fixed graph $G$. There are two players, a
set of $c$ cops, for some integer $c\ge 1$, and a robber.  Initially,
the cops are placed onto vertices of their choice in $G$ (where more
than one cop can be placed at a vertex). Then the robber, being fully
aware of the cops' placement, positions himself in one of the vertices
of $G$. Then the cops and the robber move in alternate rounds, with
the cops moving first; however, players are permitted to remain
stationary on their turn if they wish. The players use the edges of
$G$ to move from vertex to vertex. The cops win and the game ends if
eventually a cop steps into the vertex currently occupied by the
robber; otherwise, i.e., if the robber can elude the cops
indefinitely, the robber wins.

The {\em cop number}\/ of $G$, denoted by $c(G)$, is the minimum
number of cops needed to win on $G$.  This parameter was introduced by
Aigner and Fromme \cite{AF}, and there is now an extensive literature
on this fascinating problem. We direct the reader to the surveys
\cite{Alspach}, \cite{Hahn}, and \cite{FT} for detailed accounts of
the known results.  All results focus on connected graphs, because the
problem for a disconnected graph obviously decomposes into the sum of
the answers for each connected component.

The most well known open question in this area is Meyniel's
conjecture, published by Frankl in \cite{Frankl1}.  It states that for
every connected graph $G$ on $n$ vertices, $O(\sqrt{n})$ cops are
enough to win. This conjecture, if true, is best possible, as
projective plane graphs ($n$-vertex graphs without cycles of length 3
and 4, and with all degrees at least $c\sqrt{n}$) are easily seen to
require at least $c\sqrt{n}$ cops. So far, the progress towards
establishing Meyniel's conjecture has been rather slow. Frankl
\cite{Frankl1} proved the upper bound of $O(n\log\log n/\log n)$; some
twenty years later Chiniforooshan \cite{C} improved it to $O(n/\log
n)$.  Finally, the upper bound of $n / 2^{(1-o(1))\sqrt{\log n}}$ was
established by Lu and Peng \cite{LP}, and very recently Scott and
Sudakov \cite{SS} posted an alternative proof of this result.  Bounds
on the typical behavior of the cop number for the random graph $\Gnp$
have been obtained for various values of $p=p(n)$ by Bollob\'as, Kun
and Leader \cite{BKL} and by \L uczak and Pra\l at \cite{LP}. Many
other versions and ramifications of the above described classical
setting have been studied, such as the ranged version in
\cite{BCP-range}, limited visibility in \cite{IKK-limited-vis}, etc.,
but we do not pursue them here.

We employ an approach based on the notion of \emph{expansion}, which
has had many applications in mathematics and theoretical computer
science.  Recall that a graph is said to be a {\em $c$-expander} if
every subset $S$ of at most $n/2$ vertices has $|N(S) \setminus S| >
c|S|$. In particular, using this approach we are able to provide
another proof for the Lu-Peng bound mentioned above.

It also allows us to address the directed graph version of the Cops
and Robbers problem. The setting here is a straightforward adaptation
of the undirected setting described above, with the only difference
being that the players need to respect the direction of any edge while
moving along it. Problems for directed graphs (digraphs) are usually
much more difficult.  To the best of our knowledge, there have been no
results on this problem in this case.  In Section \ref{sec:digraph},
we observe that the essence of the problem is to consider only
\emph{strongly connected}\/ digraphs, i.e., those which have directed
paths from any vertex to any other vertex. We prove the following
general upper bound.

\begin{theorem}
  \label{thm:digraph}
  Every strongly connected digraph on $n$ vertices has cop number at most
  $O\big(n \cdot \frac{(\log \log n)^2}{\log n} \big)$.
\end{theorem}

We then use a purely expansion-based argument to provide an
alternate proof of Lu and Peng's result for general graphs
\cite{LP}.

\begin{theorem}
  \label{thm:general}
  Every connected graph on $n$ vertices has cop number at most
  $n / 2^{(1-o(1)) \sqrt{\log_2 n}}$.
\end{theorem}

Our approach has the added advantage that it still works in the case
when the robber moves faster than the cops.  Indeed, this setting
was recently considered by \cite{FGKNS-fast}.  (It is not an
interesting problem if a cop can move faster than the robber,
because then one cop is sufficient: he can chase down the robber.)
So we consider the case when the robber moves at speed $R > 1$ and
cop moves at speed 1; the robber can take any walk of length $R$
from his current position, but he is not allowed to pass through any
vertex occupied by a cop. With our alternate approach, we are able
to prove an result analogous to that of Lu and Peng, but for a faster
robber.

\begin{theorem}
  \label{thm:fast-upper}
  Let $R \geq 1$ be a given finite constant, and let $\alpha = 1 +
  \frac{1}{R}$.  In every connected graph on $n$ vertices, $n /
  \alpha^{(1-o(1)) \sqrt{\log_\alpha n}}$ cops are sufficient to catch
  any robber who moves at speed $R$.
\end{theorem}

\noindent \textbf{Remark.}\, Observe that for the original case $R=1$,
the constant $\alpha$ is precisely 2.  Therefore, this extends all
current best results in the traditional setting.

\vspace{3mm}

It is also interesting to note that in the fast robber setting, the
cop number can be drastically different.  Indeed, Proposition
\ref{prop:speed-2} in Section \ref{sec:fast-upper} exhibits an
$n$-vertex graph for which the cop number jumps from 2 to
$\Theta(\sqrt{n})$ when the robber's speed increases from 1 to 2.
For higher speeds, we also show that the general lower bound climbs
beyond $n^{1/2}$, and even reaches $\Omega(n)$ for an infinite-speed
robber.

\begin{theorem}
  \label{thm:fast-lower}
  For any given robber speed $R > 2$, the following hold for
  sufficiently large $n$.
  \begin{description}
  \item[(i)] If $R < \infty$, there exists an $n$-vertex graph which
    requires at least $n^{1-\frac{1}{R-2}}$ cops.
  \item[(ii)] If $R = \infty$, there exists an $n$-vertex graph which
    requires at least $\frac{n}{800^2}$ cops.
  \end{description}
\end{theorem}

Throughout our paper, we will omit floor and ceiling signs whenever
they are not essential, to improve clarity of presentation.  All
logarithms are in base $e \approx 2.718$ unless otherwise specified.
The following asymptotic notation will be utilized extensively.  For
two functions $f(n)$ and $g(n)$, we write $f(n) \ll g(n)$, $f(n) =
o(g(n))$, or $g(n) = \omega(f(n))$ if $\lim_{n \rightarrow \infty}
f(n)/g(n) = 0$, and $f(n) = O(g(n))$ or $g(n) = \Omega(f(n))$ if
there exists a constant $M$ such that $|f(n)| \leq M|g(n)|$ for all
sufficiently large $n$. The number of vertices $n$ is assumed to be
sufficiently large where necessary.

\section{Preliminaries}

Previous attempts to solve the general case of this problem have
relied on the following two observations.  Recall that $c(G)$ denotes
the cop number of $G$.

\begin{lemma}
  \label{lem:degree-diameter}
  Let $G$ be a connected $n$-vertex graph.
  \begin{description}
  \item[(i)] If $v$ is a vertex of maximum degree $\Delta$, then $c(G)
    \leq 1 + c(G')$, where $G'$ is a connected graph with at most
    $n-1-\Delta$ vertices.
  \item[(ii)] If $v_1 v_2 \ldots v_t$ is a geodesic of $G$, then $c(G)
    \leq 1 + c(G')$, where $G'$ is a connected graph with at most
    $n-t$ vertices.
  \end{description}
\end{lemma}

\noindent \textbf{Proof.}\, For part (i), permanently station a cop at
$v$.  Let $G[U_1], \ldots, G[U_k]$ be the connected components of $G
\setminus (\{v\} \cup N(v))$.  Since the robber can never enter $\{v\}
\cup N(v)$, he must stay in whichever component he starts in, say
$G[U_j]$.  The cops can pass through each other, so as long as there
are $\max_i c(G[U_i])$ other cops, they can all move into $G[U_j]$ and
capture the robber there.  As $\sum_i |U_i| = n - 1 - \Delta$, this
completes the first part.

For part (ii), a result of Aigner and Fromme \cite{AF} shows that
one cop is sufficient to patrol $\{v_1, \ldots, v_t\}$.  Specifically,
one cop can ensure that after finitely many moves, the robber is never
in $\{v_1, \ldots, v_t\}$.  Then, the same argument as above shows
that if the connected components of $G \setminus \{v_1, \ldots, v_t\}$
are $G[U_1], \ldots, G[U_k]$, then only $\max_i c(G[U_i])$ more cops
are required to capture the robber.  We will not mention their proof
here, because we do not use this part anywhere in our approach.
However, it should be noted that their proof relies critically on the
bidirectionality of the edges, and the ordinary speed of the robber.
\hfill $\Box$

\vspace{3mm}

Our new approach is to use expansion.  Let $\del S$ denote the set of
all vertices which are outside $S$, but adjacent to some vertex in
$S$.

\begin{lemma}
  \label{lem:exp}
  Let $G$ be a connected $n$-vertex graph.  Suppose that $G$ has a set
  of vertices $S$ with $|S| < n/2$, but $|\del S| \leq p|S|$ for some
  $0 < p < 1$.  Then $c(G) \leq p|S| + c(G')$, where $G'$ is a
  connected graph with at most $n - |S|$ vertices.  This holds even if
  the robber is permitted to move at speed $R > 1$.
\end{lemma}

\noindent \textbf{Proof.}\, Permanently station one cop on each vertex
in $\del S$.  Let $G[U_1], \ldots, G[U_k]$ be the connected components
of $G \setminus \del S$.  The barrier of cops will prevent the robber
from ever entering $\del S$, so in particular, he will be forced to
remain within a single connected component $G[U_j]$.  Therefore, as in
the proof of Lemma \ref{lem:degree-diameter}, only $\max_i c(G[U_i])$
more cops are required to capture the robber.

It remains to show that all $|U_i| \leq n - |S|$.  Observe that since
$|S| < n/2$, every connected component spanned by $S$ has size at most
$|S| < n/2 < n - |S|$.  On the other hand, every connected component
spanned by $G \setminus (S \cup \del S)$ obviously has size at most $n
- |S| - |\del S| \leq n - |S|$.  This completes the proof.  \hfill
$\Box$

\vspace{3mm}

Although all statements in this paper are about deterministic graphs,
we will use the probabilistic method to develop strategies and provide
constructions.  Therefore, we need a few basic probabilistic
statements.  Let us recall the Chernoff bound (see, e.g., \cite{AS}).

\begin{fact}
  \label{fact:chernoff}
  For any $\epsilon > 0$, there exists $c_\epsilon > 0$ such that any
  binomial random variable $X$ with mean $\mu$ satisfies
  \begin{displaymath}
    \pr{X < (1-\epsilon) \mu} < e^{-\frac{\epsilon^2}{2} \mu}
    \andsp
    \text{and}
    \andsp
    \pr{X > (1+\epsilon) \mu} < e^{-c_\epsilon \mu},
  \end{displaymath}
  where $c_\epsilon$ is a constant determined by $\epsilon$.  When
  $\epsilon = 1/2$, we may take $c_{1/2} = \frac{1}{10}$.
\end{fact}

\section{Directed graphs}
\label{sec:digraph}

Recall that a digraph is \emph{strongly connected}\/ if there is a
directed path from any vertex to any other vertex, \emph{weakly
  connected}\/ if its underlying undirected graph is connected, and
\emph{disconnected}\/ if its underlying undirected graph is
disconnected.  We claim that the essence of the directed case of this
problem is to investigate the cop number of an arbitrary strongly
connected digraph $D$.  Indeed, as in the case of ordinary graphs, if
the underlying undirected graph $G$ is disconnected, with connected
components $G[V_1]$, \ldots, $G[V_t]$, then the cop number of $D$ is
clearly the sum of the cop numbers of $D[V_1]$, \ldots, $D[V_t]$.  The
following proposition shows that the weakly connected case reduces to
solving the strongly connected case.

\begin{proposition}
  \label{prop:digraph-weak-strong}
  Let $D$ be a directed graph, whose strongly connected components are
  $D[V_1]$, \ldots, $D[V_t]$.  Let $c_i$ denote the cop number of
  $D[V_i]$, and construct a directed acyclic graph $D'$ with vertex
  set $[t]$, where $\overrightarrow{ij}$ is an edge if and only if $D$
  has an edge from $V_i$ to $V_j$.  Then the problem of determining
  the cop number of $D$ reduces to an optimization problem involving
  only $D'$ and the $c_i$.
\end{proposition}

\noindent \textbf{Proof.}\, In a directed acyclic graph, we call
$v$ a \emph{source vertex}\/ if it has in-degree zero, and we say
that $v$ \emph{feeds into}\/ $w$ if there is a directed path from
$v$ to $w$. The first observation is that it is never useful to
initially position cops in any $V_i$ where $i$ is a non-source
vertex of $D'$.  Indeed, consider positioning those cops in a strong
component $V_i'$, where $i'$ is  a source vertex of $D'$ which feeds
into $i$ instead. Let $V_j$ be the strong component containing the
robber's initial vertex.  If $i$ does not feed into $j$, then this
alternate placement makes no difference, because the cops in $V_i$
would be useless anyway. Otherwise, if $i$ feeds into $j$, let all
cops initially stay stationary until the relocated cops move to
their old positions in $V_i$.  Then, run the old algorithm (which
had those cops starting in $V_i$).

Therefore, we only need to choose the numbers of cops to place in
the strong components corresponding to source vertices of $D'$.
Assign an integer variable $x_i$ to each source.  Consider any
vertex $j \in D'$, and let $S$ be the set of source vertices which
feed into it.  We must have the inequality $\sum_S x_i \geq c_j$,
because if the robber started in $V_j$, then the only sources of
cops are from the $V_i$ with $i \in S$.  We thus introduce one
constraint per vertex of $D'$.

It remains to show that if all constraints are satisfied, then the
robber certainly can be caught.  Let the robber's initial position be
in $V_j$.  Route all possible cops into $V_j$.  By the constraints, we
will be able to move at least $c_j$ cops there.  If the robber stays
in $V_j$, he will certainly be caught eventually, so he must exit to
another strong component $V_k$.  In $D'$, the set of sources which
feed into $k$ is a superset of those which feed into $j$, so we can
move the cops in $V_j$ to $V_k$, and move the cops from the other
sources directly to $V_k$.  This will force the robber out of $V_k$,
and the process must terminate eventually because $D'$ is finite.
\hfill $\Box$

\vspace{3mm}

Solving the optimization problem is outside of the scope of this
paper, since it has an entirely different flavor.  Instead, we now
proceed to prove Theorem \ref{thm:digraph}.  Diameter-based
arguments break down completely, so previous bounds for general
graphs (e.g., Frankl \cite{Frankl1}, Chiniforooshan \cite{C}, and Lu
and Peng \cite{LP}) do not apply.  However, expansion is immune to
this difficulty.  In the context of directed graphs, let us call a
digraph a {\em $c$-in-expander} if every subset $S$ of at most $n/2$
vertices has $|\del^- S| \geq c|S|$, where $\del^- S$ is the set of
all vertices $v \not \in S$ which have a directed edge into $S$. Let
us bring some basic tools from the previous section to the directed
setting.

\begin{lemma}
  \label{lem:di-degree-exp}
  Let $D$ be a strongly connected $n$-vertex digraph.
  \begin{description}
  \item[(i)] If $v$ is a vertex of out-degree $\Delta$, then $c(D)
    \leq 1 + c(D')$, where $D'$ is a strongly connected digraph with
    at most $n-1-\Delta$ vertices.
  \item[(ii)] Suppose that $D$ has a set of vertices $S$ with $|S| <
    n/2$, but $|\del^- S| \leq p|S|$ for some $0 < p < 1$.  Then $c(D)
    \leq p|S| + c(D')$, where $D'$ is a strongly connected digraph
    with at most $n - |S|$ vertices.
  \end{description}
\end{lemma}

\noindent \textbf{Proof.}\, Part (i) has essentially the same proof as
Lemma \ref{lem:degree-diameter}(i).  By permanently stationing one cop
on $v$, the robber is never able to enter $\{v\} \cup N^+(v)$.  Let
$D[U_1], \ldots, D[U_k]$ be the strongly connected components of $D
\setminus (\{v\} \cup N^+(v))$.  Since these are strong components,
observe that if the robber ever moves out of some $D[U_i]$, then it
can never return to $D[U_i]$.  Therefore, with only $\max_i c(D[U_i])$
additional cops, the robber will eventually be chased out of all
strong components until he is trapped and captured in the final one.
The proof of part (ii) is analogous, and follows the same argument as
Lemma \ref{lem:exp}.  \hfill $\Box$

\vspace{3mm}

The previous lemma is cumbersome to apply by itself.  However, it
allows us to clean up our graph, at the cost of reserving a few cops
for this purpose.  We record the following statement, which is more
convenient to use.

\begin{corollary}
  \label{cor:digraph-nice}
  Let $D$ be a strongly connected digraph with $n$ vertices, and let
  $0 < p < 1$ be arbitrary.  Then $c(D) \leq pn + c(D')$, where $D'$
  is a strongly connected digraph with at most $n$ vertices and
  maximum out-degree at most $1/p$, which is also a $p$-in-expander.
\end{corollary}

\noindent \textbf{Proof.}\, We repeatedly apply Lemma
\ref{lem:di-degree-exp}.  As long as there is a vertex $v$ of
out-degree at least $\frac{1}{p}$, part (i) shows that at the cost of
one cop, we can reduce the number of vertices by at least $\frac{1}{p}
+ 1$.  Similarly, if there is a set $S$ of at most half the vertices
with $|\del^- S| \leq p|S|$, we can reduce the number of vertices by
at least $|S|$, at the cost of $p|S|$ cops.  Note that in both cases,
the number of cops expended is at most a $p$-fraction of the number of
vertices discarded.  Therefore, if we repeat this process until
exhaustion, we will have a digraph $D'$ with $m \leq n$ vertices, with
the stated properties, such that $c(D) \leq p(n-m) + c(D') \leq pn +
c(D')$, as claimed. \hfill $\Box$

\vspace{3mm}

We are now ready to prove Theorem \ref{thm:digraph}.  Corollary
\ref{cor:digraph-nice} implies that it is an immediate consequence of
the following final lemma.

\begin{lemma}
  Let $p = \frac{13 (\log \log n)^2}{\log n}$.  Every strongly
  connected digraph $D$ on $m \leq n$ vertices with maximum out-degree
  at most $1/p$ and in-expansion at least $p$ can be guarded by at
  most $2pn$ cops.
\end{lemma}

\noindent \textbf{Proof.}\, Note that if $m \leq pn$, we are trivially
done by placing a cop on each vertex of $D$.  So, we may assume that
$m > pn$.  Let $r = \frac{6}{p} \log \frac{4}{p}$.  For each vertex
$v$, let $B^+(v)$ denote the set of all vertices which are reachable
from $v$ by (directed) walks of length at most $r$.  Similarly, for $S
\subset V$, let $B^-(S)$ contain every vertex that can reach some
vertex in $S$ by a directed walk of length at most $r$.

Our first claim is that it is possible to position $2pn$ cops so
that for every subset $S$ of size $1 \leq |S| \leq 2p^{-r}$, the set
$B^-(S)$ contains at least $|S|$ cops.  Indeed, Inequality
\ref{ineq:di-p} from the Appendix gives
\begin{displaymath}
  2p^{-r} (1+p)^r
  \leq
  np/2
  <
  m/2,
\end{displaymath}
so for such a set $S$, the in-expansion property ensures that
$|B^-(S)| \geq |S|(1+p)^r$.  Note that Inequality \ref{ineq:di-p} from
the Appendix also shows that $(1+p)^r \geq \frac{16}{p} \log n$.

Therefore, if we position cops randomly, by independently placing a
cop at each vertex with probability $p$, the expected number of cops
in $B^-(S)$ is at least $|S| \cdot 16 \log n$.  The Chernoff bound
(Fact \ref{fact:chernoff}) shows that the probability that this is
below half its expectation is at most $e^{-\frac{1}{8} |S|\cdot 16
\log n} \leq n^{-2 |S|}$.

Since the number of subsets of $s$ vertices is at most $n^s$, a
union bound over all $S$ of size $s$ shows that with probability at
least $1 - n^{-s}$, every such $B^-(S)$ contains at least $|S|\cdot
8 \log n \geq |S|$ cops.  Taking another union bound over all $s \in
\{1, \ldots, 2p^{-r}\}$, we see that \whp, this holds for every $1
\leq |S| \leq 2p^{-r}$.  Also, the Chernoff bound implies that \whp,
at most $2pn$ cops were placed by the random process.  Putting these
two together, we see that it is indeed possible to place only $2pn$
cops so that for every $S$ of size $1 \leq |S| \leq 2p^{-r}$, the
set $B^-(S)$ contains at least $|S|$ cops.

Now assume that the cops are placed as above.  Let the robber's
position be $v$.  By the maximum out-degree condition, $|B^+(v)|
\leq 1+p^{-1}+\ldots + p^{-r}< 2p^{-r}$.  We use Hall's theorem to
show that for each $w \in B^+(v)$, there is a distinct cop $c_w$
which can reach it within $r$ moves.  The necessary condition is
precisely what we established by the above argument.  Therefore, by
sending each $c_w$ to position $w$, the robber will definitely be
captured within $r$ moves. \hfill $\Box$

\section{General graphs}

In the previous section, we used Hall's theorem to route distinct cops
to each position which needed to be blocked.  As in the paper of Lu
and Peng \cite{LP}, this argument can be improved by performing
several iterations.  The main idea is to draw conclusions from the
failure of the Hall condition.  The $k=1$ version of the following
lemma essentially appears in \cite{LP}, but our statement has an
expansion flavor built in, and so is more amenable to our approach.

\begin{lemma}
  \label{lem:split-hall}
  Given a bipartite graph with parts $A$ and $B$ and a number $k$, it
  is always possible to partition $A = S \cup T$ such that $|N(S)|
  \leq k|S|$, and for every subset $U \in T$, we have $|N(U)| \geq
  k|U|$.
\end{lemma}

\noindent \textbf{Proof.}\, Start with $S = \emptyset$.  As long as
there is a subset $U \subset A$ with $|N(U)| \leq k|U|$, add $U$ to
$S$, and delete $U$ from $A$.  It is clear that at the end of this
process, if we consider the original graph, $|N(S)| \leq k|S|$.
However, in our modified graph, every subset of the remaining $A$
expands by at least $k$ times in $B$.  \hfill $\Box$

\vspace{3mm}

Next, we isolate a component of our directed graph proof, so that we
can use it in a modular form.  Let $B_r(S)$ denote the set of all
vertices which are within distance $r$ from at least one vertex in
$S$.

\begin{lemma}
  \label{lem:sprinkle}
  Let $n, p, r$ be given, with $np$ sufficiently large.  In every
  $n$-vertex graph $G$, it is possible to distribute $2pn$ cops such
  that for every set $S$ with $|B_r(S)| \geq \frac{16}{p} |S| \log n$,
  there are at least $|S|$ cops in $B_r(S)$.
\end{lemma}

\noindent \textbf{Proof.}\, Position cops randomly, by independently
placing a cop at each vertex with probability $p$.  For each $S$ in
the statement, the expected number of cops in $B_r(S)$ is at least
$16 |S| \log n$.  The Chernoff bound (Fact \ref{fact:chernoff})
shows that the probability that this is below half its expectation is
$e^{-\frac{1}{8} \cdot 16 |S| \log n} \leq n^{-2 |S|}$.

Since the number of subsets of $s$ vertices is at most $n^s$, a
union bound over all $S$ of size $s$ shows that with probability at
least $1 - n^{-s}$, every such $B_r(S)$ contains at least $8 |S|
\log n \geq |S|$ cops.  Taking another union bound over all $s \in
\{1, \ldots, n\}$, we see that with probability at least
$1-\frac{2}{n}$, this holds for every $1 \leq |S| \leq n$.  Yet
$\bin{n, p}$ is at most $2np$ \whp\ by the
Chernoff bound, so we conclude that there is positive probability of
our procedure giving all of the desired properties, using only $2pn$
cops.  \hfill $\Box$

\vspace{3mm}

Next, we translate Corollary \ref{cor:digraph-nice} to the case of
undirected graphs, via Lemmas \ref{lem:degree-diameter}(i) and
\ref{lem:exp}.  The proof of the following statement is analogous to
Corollary \ref{cor:digraph-nice}, so we do not record it again.

\begin{corollary}
  \label{cor:nice}
  Let $G$ be a connected graph with $n$ vertices, and let $0 < p < 1$
  be arbitrary.  Then $c(G) \leq pn + c(G')$, where $G'$ is a
  connected graph with $m \leq n$ vertices and maximum degree at most
  $1/p$, which is also a $p$-expander.
\end{corollary}

In light of this corollary, it is clear that Theorem \ref{thm:general}
is immediate from the following final lemma.

\begin{lemma}
  \label{lem:general}
  There is a function $p = p(n) = 2^{-(1-o(1)) \sqrt{\log_2 n}}$ for
  which the following holds.  Every connected graph $G$ on $m \leq n$
  vertices with maximum degree less than $1/p$ and expansion at least
  $p$ can be guarded by at most $(1+o(1)) \sqrt{\log_2 n} \cdot 2pn$
  cops.
\end{lemma}

\noindent \textbf{Proof.}\, As in the proof for directed graphs, we
may assume that $m > pn$, or else we are trivially done.  Inequality
\ref{ineq:p} shows that there is a function $p = 2^{-(1-o(1))
  \sqrt{\log_2 n}}$ and a positive integer $l = (1+o(1)) \sqrt{\log_2
  n}$ such that when we define $k = \frac{16}{p} \log n$, we have the
inequalities
\begin{equation}
  \label{ineq:p-repeat}
  k^{l+1} \leq (1+p)^{-2^l} np/2
  \andsp
  \text{and}
  \andsp
  (1+p)^{2^l} \geq k.
\end{equation}
We will split the cops into $l+1$ groups $C_0, C_1, \ldots, C_l$, each
of size $2pn$.  Choose the initial positions of the cops in $C_i$ by
applying Lemma \ref{lem:sprinkle} with parameter $r = 2^i$.  Let the
robber's initial position be $v$.

Let $N_0 = B_1(v)$ be the set of vertices that the robber can reach in
1 move.  By the maximum degree condition, $|N_0| \leq \frac{1}{p} <
k$.  Consider the auxiliary bipartite graph in which $A = N_0$, $B =
V$, and $a$ is adjacent to $b$ if and only if they are at distance at
most 1 in $G$.  Then Lemma \ref{lem:split-hall} implies that we can
partition $N_0 = S_0 \cup T_0$ such that (in $G$) $|B_1(S_0)| \leq k
|S_0|$ and every subset $U \subset T_0$ has $|B_1(U)| \geq k |U|$.
Therefore, by construction of $C_0$, Hall's theorem shows us how to
send a distinct cop from $C_0$ to each vertex of $T_0$ in the first
move, preventing the robber from ever occupying a vertex of $T_0$.

Thus the robber's position after his second move is restricted to
$N_1 = B_1(S_0)$, and $|B_1(S_0)|
  \leq
  k |S_0|
  \leq
  k^2$.
Repeating the same trick, we can partition $N_1 = S_1 \cup T_1$ and
use the cops in $C_1$ to prevent the robber from ever entering $T_1$,
yet $|B_2(S_1)| \leq k |S_1| \leq k^3$.  Hence the
robber's position after his 4th move is restricted to $N_2 =
B_2(S_1)$.

The radii of the balls double at each iteration of this procedure, so
we eventually conclude that after his $2^l$-th move, the robber is
still contained within a set $N_l = B_{2^{l-1}}(S_{l-1})$ of size at
most $k^{l+1}$.  However, when we iterate the argument a final time,
the partition $N_l = S_l \cup T_l$ must have $S_l = \emptyset$.
Indeed, since $G$ is a $p$-expander, every non-empty set $S$ of size
at most $(1+p)^{-(r-1)} m/2 > (1+p)^{-(r-1)} np/2$ has $|B_r(S)| \geq
(1+p)^r |S|$.  As Inequality \eqref{ineq:p-repeat} ensures that $|N_l|
\leq k^{l+1} \leq (1+p)^{-2^l} np/2$ and $(1+p)^{2^l} \geq k$, we
conclude that $S_l$ is indeed empty.  Therefore, the cops in $C_l$ can
completely cover $N_l$ within $2^l$ moves.  Since $N_l$ was the set of
possible positions for the robber after his $2^l$-th move, the robber
is captured.  \hfill $\Box$

\section{Fast robber}
\label{sec:fast-upper}

In this section, we assume that the robber can traverse up to $R$
edges in a single move.  Cops may only move by a single edge per move.
We begin by observing that the cop number of a graph can dramatically
increase even if the robber's speed only grows to $R=2$.

\begin{proposition}
  \label{prop:speed-2}
  Let $G$ be the 1-subdivision of $K_n$, where a vertex is added on
  each edge.  The ordinary cop number of $G$ is 2, but if the robber
  can move at speed 2, then the cop number rises to $\lceil n/2 \rceil
  = \Theta(\sqrt{v(G)})$.
\end{proposition}

\noindent \textbf{Proof.}\, Call a vertex of a 1-subdivision an
\emph{internal vertex}\/ if it was added to subdivide an edge, and a
\emph{join vertex}\/ otherwise.  In the ordinary setting, by placing
two cops on arbitrary join vertices $a, b$, they can catch the robber
within 3 moves.  Indeed, if the robber starts on the internal vertex
between two join vertices $u, w$, then both cops move toward $u, w$ on
their first move.  Regardless of which of $u, w$ the robber moves to,
a cop will be adjacent, and can catch him on the next turn.
Otherwise, if the robber starts on a join vertex $v$, then the cop at
$a$ moves to the internal vertex between $a$ and $v$.  The robber must
move to an internal vertex, say between $v$ and $w$.  The cop at $a$
follows him to $v$, and the other cop moves to the internal vertex
between $b$ and $w$.  The robber will now be caught in the next round.

On the other hand, if the robber moves at speed 2, then $\lceil n/2
\rceil$ cops are required to catch him.  To see this, note that any $m
< \lceil n/2 \rceil$ cops can be immediately adjacent to only at most
$2m < n$ join vertices.  So, the robber can choose a non-dominated
join vertex, say $v$, to start on and wait.  When a cop moves adjacent
to him, there will be at least one join vertex $w$ with no adjacent
cop.  Importantly, the vertex between $v$ and $w$ is unoccupied, since
otherwise a cop would be adjacent to $w$.  So, the robber can advance
to $w$ in a single move, and again be nonadjacent to any cop.  He can
repeat this indefinitely, eluding $\lceil n/2 \rceil - 1$ cops.

Note that if $\lceil n/2 \rceil$ cops are used, they can initially sit
on internal vertices so that all join vertices are dominated.  The
robber must then select an internal vertex for his initial position,
say between join vertices $v$ and $w$.  These two vertices are not
dominated by the same cop, because the only vertex which does is
occupied by the robber.  So, the two cops which dominate $v$ and $w$
can advance to occupy $v$ and $w$ in their first turn.  This traps the
robber, and he will be captured in the next round.  \hfill $\Box$

\vspace{3mm}

Let us now turn our attention to upper bounds.  Unfortunately,
diameter-based arguments completely break down, because Lemma
\ref{lem:degree-diameter}(ii) does not hold for fast robbers.
However, since even a fast robber cannot pass through vertices
occupied by cops, Lemma \ref{lem:exp} still applies.  Therefore, we
can adapt the proof from the previous section to this case.  The first
step is to extend Lemma \ref{lem:degree-diameter}(i) to this setting.

\begin{lemma}
  \label{lem:fast-degree}
  Let $n, p, R$ be given, with $np$ sufficiently large.  Every
  $n$-vertex graph $G$ has a set $U$ of $2pn$ vertices such that the
  following holds.  Place $R$ cops on each vertex of $U$, and let the
  robber choose a starting position.  Then there is a set $S$ of size
  at most $\big( \frac{16}{p} \log n \big)^{2^R}$ such that (since the
  cops move first) the robber's position after his first move must lie
  in $S$.
\end{lemma}

\noindent \textbf{Proof.}\, Let $k = \frac{16}{p} \log n$.  We
construct $U$ such that that every vertex $v$ with $|B_1(v)| \geq k$
has at least $|B_1(v)| \cdot \frac{p}{2}$ vertices of $U$ in $N(v)$.
By independently including each vertex with probability $p$, the
probability that this property fails at a fixed $v$ is at most
$e^{-kp/8} = n^{-2}$ by the Chernoff bound.  Combining a union bound
over all $v$ with the fact that $\bin{n, p}$ is at most $2np$ \whp, we
see that we have positive probability of obtaining the desired
construction.

Now let $C_1, \ldots, C_R$ be $R$ sets of cops, where each set has one
cop on each vertex of $U$.  The robber cannot select an initial vertex
with degree at least $k$, or else he will be adjacent to a cop in
$C_1$ (who will catch him immediately, since cops move first).  So,
assume that the robber's initial vertex $v$ has $|B_1(v)| \leq k$.

We will simultaneously dispatch the cops in $C_2, \ldots, C_R$, so
that in their first move, they occupy high-degree vertices in the
vicinity of $v$.  Consider the vertices of $B_1(v)$ which have degree
at least $\frac{2}{p} k$.  By construction, each such vertex will have
at least $k \geq |B_1(v)|$ cops in $C_2$ in its neighborhood, so by
the greedy algorithm, we may send cops in $C_2$ to occupy these
vertices before the robber has a chance to move.  Since the robber
cannot pass through any cops, he must avoid these vertices forever.
Let $S_1 \subset B_1(v)$ be the remaining vertices, and let $S_2 =
B_1(S_1)$.  These are the potential positions that the robber can
reach within distance 2, and we have $|S_2| \leq \frac{2}{p} k^2$.

Repeating this argument again, we see that we may send cops in $C_3$
to occupy vertices in $S_2$ that have degree at least $\frac{2}{p}
\cdot \frac{2}{p} k^2 = \big( \frac{2}{p} \big)^2 k^2$.  Then, $S_3 =
B_1(S_2)$ is the set of potential positions that the robber can reach
within distance 3, and $|S_3| \leq \big( \frac{2}{p} \big)^2 k^2 \cdot
|S_2| \leq \big( \frac{2}{p} \big)^3 k^4$.  Continuing in this way, we
see that after dispatching $C_R$, we have restricted the set of
positions that the robber can reach within distance $R$ to a set of
size at most
\begin{displaymath}
  \left( \frac{2}{p} \right)^{2^{R-1}-1} k^{2^{R-1}} < k^{2^R}.
\end{displaymath}
This is the maximum distance he can cover in his first move.  \hfill
$\Box$

\vspace{3mm}

Now we are ready to extend our earlier proof to the fast robber
setting.  Since Lemma \ref{lem:exp} holds for fast robbers, the
obvious translation of Corollary \ref{cor:digraph-nice} implies that
Theorem \ref{thm:fast-upper} is a consequence of the following lemma.

\begin{lemma}
  Let $R$ be a positive integer, and let $\alpha = 1 + \frac{1}{R}$.
  There is a function $p = p(n) = \alpha^{-(1-o(1)) \sqrt{\log_\alpha
      n}}$ for which the following holds.  In every connected graph
  $G$ on $m \leq n$ vertices with maximum degree less than $1/p$ and
  expansion at least $p$, $(1+o(1)) \sqrt{\log_\alpha n} \cdot 2pn$
  cops can always capture a speed-$R$ robber.
\end{lemma}

\noindent \textbf{Proof.}\, As usual, we may assume that $m > pn$, or
else we are trivially done.  Define the sequence $d_0, d_1, \ldots$
via the recursion $d_0 = R$, $d_{i+1} = d_i + \big\lceil \frac{d_i}{R}
\big\rceil$.  Let $r_i = \big\lceil \frac{d_i}{R} \big\rceil$.
Inequality \ref{ineq:fast-p} shows that there is a function $p =
\alpha^{-(1-o(1)) \sqrt{\log_\alpha n}}$ and a positive integer $l =
(1+o(1)) \sqrt{\log_\alpha n}$ such that when we define $k =
\frac{16}{p} \log n$, we have the inequalities
\begin{equation}
  \label{ineq:fast-p-repeat}
  k^{2^R} \cdot k^l \leq (1+p)^{-r_l} \frac{np}{2}
  \andsp
  \text{and}
  \andsp
  (1+p)^{r_l} \geq k.
\end{equation}
We use Lemma \ref{lem:fast-degree} to distribute $2Rpn$ cops such that
the robber's position after his first move will always be contained in
a set $N_0$ of size at most $k^{2^R}$.  We split the remaining cops
into $l+1$ groups $C_0, C_1, \ldots, C_l$.  Choose the initial
positions of the cops in $C_i$ by applying Lemma \ref{lem:sprinkle}
with parameter $r_i$.  Let the robber's initial position be $v$.

The rest of the proof is nearly identical to that of Lemma
\ref{lem:general}.  At each step, we consider the robber's set of
possible intermediate positions in $B_{d_i}(v)$, which he occupies on
or after his $\lfloor d_i / R \rfloor$-th move, but strictly before
the completion of his $(\lfloor d_i / R \rfloor + 1)$-st move.  We let
this set be $N_i$, and inductively assume it has size at most $k^{2^R}
\cdot k^i$.

Since the cops move first, they can travel by distance $\lfloor d_i /
R \rfloor + 1 \geq r_i$ by this time.  Lemma \ref{lem:split-hall}
partitions $N_i = S_i \cup T_i$ such that $|B_{r_i}(S_i)| \leq k
|S_i|$ and every subset $U \subset T_i$ has $|B_{r_i}(U)| \geq k |U|$.
By construction of $C_i$, Hall's theorem shows us how to send a
distinct cop from $C_i$ to each vertex of $T_i$.  Hence the robber
actually cannot occupy $T_i$ anytime after his $\lfloor d_i / R
\rfloor$-th move.  Therefore, the set of positions in $B_{d_{i+1}}(v)$
which the robber may occupy between his $\lfloor d_{i+1} / R
\rfloor$-th and $(\lfloor d_{i+1} / R \rfloor + 1)$-st moves is
restricted to some $N_{i+1}$, of size at most $k^{2^R} \cdot k^{i+1}$.

This procedure terminates because when we partition $N_l = S_l \cup
T_l$, the expansion property ensures that $S_l = \emptyset$.  Indeed,
the inequalities in \eqref{ineq:fast-p-repeat} are precisely what are
required to show that the sets are small enough to expand, and that
their radii are large enough for the expansion factor to exceed $k$.
Therefore, the cops in $C_l$ can completely cover $N_l$ within $r_l$
moves.  Since $N_l$ was the set of possible positions for the robber
within distance $d_l$, the robber is captured.  \hfill $\Box$

\section{Lower bound for infinitely fast robber}

We now proceed to prove that the $\Omega(\sqrt{n})$ general lower
bound can be sharpened considerably in the setting when the robber
moves faster than the cops.  As a warm-up, we start with the second
part of Theorem \ref{thm:fast-lower}, which states that there are
$n$-vertex graphs on which an infinitely fast robber can always evade
$cn$ cops, for an absolute constant $c$.  The graphs will be instances
of $\Gnp$ with $p = 200/n$.  We will need some routine lemmas about
$\Gnp$.

\begin{lemma}
  \label{lem:inf-few-edges}
  Let $p = 200/n$.  Then, \whp\ every set of $s \leq 0.6n$ vertices in
  $\Gnp$ has average degree at most $0.9np$.
\end{lemma}

\noindent \textbf{Proof.}\, Note that the average degree condition is
equivalent to enforcing that each such set spans at most $\frac{s}{2}
\cdot 0.9np$ edges.  We will take a union bound, but we split the
range for $s$ into three parts.  First, consider any set $S$ of $0.3n
< s \leq 0.6n$ vertices.  The number of edges in $S$ is $\bin{{s
    \choose 2}, p}$, so $\frac{s}{2} \cdot 0.9np$ exceeds its mean by
a factor of at least 50\%.  Therefore, the Chernoff bound (Fact
\ref{fact:chernoff}) implies that the probability that $S$ fails is at
most $e^{-0.1 \cdot \frac{s^2}{2} \cdot p} \leq e^{-3s}$.  Taking a
union bound over all such $S$ with $0.3n < s < 0.6n$, we accumulate a
failure probability of at most
\begin{equation}
  \label{ineq:vis-high}
  \sum_{s = 0.3n}^{0.6n} {n \choose s} e^{-3s}
  \leq
  \sum_{s = 0.3n}^{0.6n} \left( \frac{en}{s} e^{-3} \right)^s
  \leq
  \sum_{s = 0.3n}^{0.6n} \left( \frac{e}{0.3} e^{-3} \right)^s
  \leq
  \sum_{s = 0.3n}^{0.6n} 2^{-s}.
\end{equation}

Next, we consider $s$ in the range $\log n < s \leq 0.3n$.  Here, we
use a simpler bound for the probability that a given set of $s$
vertices spans more than $\frac{s}{2} \cdot 0.9np$ edges.  Combining
this with a union bound over all $S$ of these sizes, we bound the
total failure probability in this range by
\begin{equation}
  \label{ineq:vis-mid-presub}
  \sum_{s = \log n}^{0.3n} {n \choose s} {\frac{s^2}{2} \choose \frac{s}{2} \cdot 0.9np} p^{\frac{s}{2} \cdot 0.9np}
  \leq
  \sum_{s = \log n}^{0.3n} \left(\frac{en}{s}\right)^s \left( \frac{es}{0.9np} \right)^{\frac{s}{2} \cdot 0.9np} p^{\frac{s}{2} \cdot 0.9np}
  \leq
  \sum_{s = \log n}^{0.3n} \left[
    \frac{en}{s} \cdot \left( \frac{es}{0.9n} \right)^{90}
  \right]^s.
\end{equation}
Since the quantity in the square brackets increases with $s$ (its
exponent is $+89$), we may replace $s$ with its maximum value in this
range, to obtain an upper bound of:
\begin{equation}
  \label{ineq:vis-mid}
  \sum_{s = \log n}^{0.3n} \left[
    \frac{en}{0.3n} \cdot \left( \frac{e \cdot 0.3n}{0.9n} \right)^{90}
  \right]^s
  =
  \sum_{s = \log n}^{0.3n} \left[
    \frac{e}{0.3} \cdot \left( \frac{e}{3} \right)^{90}
  \right]^s
  <
  \sum_{s = \log n}^{0.3n} 2^{-s}.
\end{equation}

For the final range $1 \leq s \leq \log n$, we may substitute $\log n$
for $s$ into inequality \eqref{ineq:vis-mid-presub}, for the same
reason as above.  So, the total failure probability in that range is
at most
\begin{equation}
  \label{ineq:vis-low}
  \sum_{s = 1}^{\log n} \left[
    \frac{en}{\log n} \cdot \left( \frac{e \log n}{0.9n} \right)^{90}
  \right]^s
  \leq
  \log n \cdot
  \frac{en}{\log n} \cdot \left( \frac{e \log n}{0.9n} \right)^{90}
  =
  O\left( \frac{\log^{90} n}{n^{89}} \right).
\end{equation}

Combining inequalities \eqref{ineq:vis-high}, \eqref{ineq:vis-mid},
and \eqref{ineq:vis-low}, we obtain the desired result.  \hfill $\Box$

\vspace{3mm}

Our next lemma allows us to delete small (but linear-size) vertex
subsets without destroying too many edges.  Let us say that an edge is
\emph{covered}\/ by a vertex subset if one of its endpoints is in the
subset.

\begin{lemma}
  \label{lem:vis-cover}
  Let $\lambda, c$ be positive real constants such that $c > e \cdot
  (e/4)^{4\lambda}$.  Then, for $p = \lambda/n$, \whp\ every set of
  $cn$ vertices in $\Gnp$ covers at most $4np \cdot cn$ edges.
\end{lemma}

\noindent \textbf{Proof.}\, Any given set of $cn$ vertices is
potentially incident to ${cn \choose 2} + c(1-c)n^2 \leq cn^2$ edges,
each of which is independently present with probability $p$.  So, the
probability that over $4np \cdot cn$ appear is at most
\begin{displaymath}
  {cn^2 \choose 4np \cdot cn} p^{4np \cdot cn}
  \leq
  \left( \frac{e}{4p} \right)^{4np \cdot cn} p^{4np \cdot cn}
  =
  \left( \frac{e}{4} \right)^{4\lambda \cdot cn}.
\end{displaymath}
Therefore, taking a union bound over all subsets of size $cn$, the
total failure probability is at most
\begin{displaymath}
  {n \choose cn} \cdot \left( \frac{e}{4} \right)^{4\lambda \cdot cn}
  \leq
  \left[ \frac{e}{c} \cdot \left( \frac{e}{4} \right)^{4\lambda} \right]^{cn}
  =
  \alpha^{cn},
\end{displaymath}
for some constant $\alpha < 1$.  Hence this probability tends to zero,
as claimed.  \hfill $\Box$

\vspace{3mm}

\noindent \textbf{Proof of Theorem \ref{thm:fast-lower}(ii).}\,
Let $c = 800^{-2}$, and let $G$ be an instance of $\Gnp$ with $p =
200/n$.  The previous lemmas show that for sufficiently large $n$, we
can ensure that $G$ has the following properties:
\begin{description}
\item[(i)] $G$ has at least $99n$ edges (Chernoff).
\item[(ii)] Every subset of $cn$ vertices covers at most $800 cn$
  edges.
\item[(iii)] Every subset of $800 cn$ vertices covers at most $800^2
  cn$ edges.
\item[(iv)] Every subset of at most $0.6n$ vertices has average
  degree at most $0.9np$.
\end{description}
We claim that these properties are enough to allow the robber to
escape $cn$ cops indefinitely.

Indeed, suppose that there are only $cn$ cops, and let $C$ be the
set of vertices that they initially occupy.  Let $C^+$ be the union
of $C$ with all immediate neighbors of vertices in $C$, and let $U$
be the complement of $C^+$.  Property (ii) shows that $|C^+| \leq
800 cn$, so by (iii), the total number of edges covered by
$C^+$ is always at most $800^2 cn = n$.  So, $G[U]$
induces at least $98n$ edges, and hence has average degree at least
$0.98 np$. Some connected component of $G[U]$ must have at least
that average degree, and (iv) shows that it must then have size at
least $0.6n$. Therefore, $G[U]$ always has a connected component of
at least this size.

The robber's strategy is to initially place himself in an arbitrary
vertex $v$ of the largest connected component of $G[U]$, which has
size at least $0.6n$.  After the cops move, let $U'$ be the
complement of the new $C^+$.  There must still be a connected
component of size at least $0.6n$ in $G[U']$; the robber selects an
arbitrary vertex $x$ in it.  Since these two large components both
have size at least $0.6n$, they must overlap in some vertex $w$.
Therefore, there is a path $P_1$ from $v$ to $w$ entirely contained
in $U$, and a path $P_2$ from $w$ to $x$ entirely contained in $U'$.
Yet even though the cops have moved, by definition of $U$, their
current positions are still outside of $U$, since $U$ excluded their
old immediate neighborhoods. Therefore, both paths $P_i$ completely
avoid all cops, so the robber can indeed move to $x$ in his turn.
This preserves the condition that he is always in the largest
connected component outside $C^+$, so he can repeat this
indefinitely.  \hfill $\Box$

\medskip

\noindent \textbf{Remark.}\, A more careful implementation of the
above argument allows the robber to escape when his speed $R$ is not
infinite, but rather at least $C\log n$ for some large enough
constant $C>0$. This is due to the fact that the large connected
subgraph of $G[U]$ in the above argument can be chosen in addition
to be of logarithmic diameter, allowing the robber to escape from it
to $U'$ in a logarithmic number of steps. A similar argument is
presented in more detail in the next section.

\section{Lower bounds for finite-speed robber}

We now extend the ideas of the previous section to prove the first
part of Theorem \ref{thm:fast-lower}.  In the last section,
connectivity alone was enough, since the robber could move
infinitely quickly.  Here, we also need to control the lengths of
the paths involved.  The graph will still be an instance of $\Gnp$,
but this time $p$ will be of order $n^{\frac{1}{R-2}-1}$.  As usual,
we begin by stating some routine facts about $\Gnp$.

\begin{lemma}
  \label{lem:fast-robber-well-connected}
  In $\Gnp$, \whp\ there is an edge between every pair of disjoint
  sets of size $s_0 = \frac{3}{p} \log n$.
\end{lemma}

\noindent \textbf{Proof.}\, For any fixed pair of disjoint sets of
size $s_0$, the probability that all crossing edges are absent is at
most $(1-p)^{s_0^2}$.  There are at most ${n \choose s_0}^2$ ways to
choose these sets, so a union bound implies that the probability that
this property does not hold in $\Gnp$ is at most
\begin{displaymath}
  {n \choose s_0}^2 (1-p)^{s_0^2}
  \leq
  \left( \frac{en}{s_0} \right)^{2s_0} e^{-p s_0^2}
  =
  \left[ \left( \frac{en}{s_0} \right)^2 e^{-p s_0} \right]^{s_0}
  \leq
  \left[ n^2 e^{-3 \log n} \right]^{s_0}
  = n^{-s_0}
  = o(1).
\end{displaymath}
\hfill $\Box$

\vspace{3mm}

\begin{lemma}
  \label{lem:fast-robbers-few-edges}
  In $\Gnp$, \whp\ every set of size $s \leq s_0 = \frac{3}{p} \log n$
  spans at most $s \cdot 6 \log n$ edges.
\end{lemma}

\noindent \textbf{Proof.}\, For fixed $s$, the number of sets of $s$
vertices is ${n \choose s}$.  The probability that a particular set of
$s$ vertices spans at least $k = s \cdot 6 \log n$ edges is at most
${s^2/2 \choose k} p^k$.  Therefore, the probability that our property
fails for a certain fixed $s$ is at most:
\begin{eqnarray*}
  {n \choose s} \cdot {s^2/2 \choose s \cdot 6 \log n} p^{s \cdot 6 \log n}
  &\leq&
  \left( \frac{en}{s} \right)^s \cdot
  \left( \frac{esp}{12 \log n} \right)^{s \cdot 6 \log n} \\
  &=&
  \left[
    \frac{en}{s} \cdot \left( \frac{esp}{12 \log n} \right)^{6 \log n}
  \right]^s \\
  &\leq&
  \left[
    en \cdot \left( \frac{e}{4} \right)^{6 \log n}
  \right]^s.
\end{eqnarray*}
Since $(e/4)^6 \approx e^{-2.3}$, this probability is at most $n^{-s}$
for large $n$.  Taking a final union bound over all $s \leq s_0$, we
see that the total failure probability is still $o(1)$, as claimed.
\hfill $\Box$

\vspace{3mm}

\begin{lemma}
  \label{lem:fast-robber-degrees-spread}
  Let $\gamma > 0$ be fixed, and suppose $np \rightarrow \infty$.
  Then $\Gnp$ has the following property \whp.  For every integer $t$
  between $\gamma np$ and $\frac{\gamma^3}{2e^5} n$, every subset
  $U$ of $t$ vertices has the property that the number of
  vertices $v\notin U$ with $d_U(v) \geq \gamma np$ is at most $3 \cdot
  \frac{t}{\gamma np}$.
\end{lemma}

\noindent \textbf{Proof.}\, Let $k = 3 \cdot \frac{t}{\gamma np}$.
For each fixed $t$, there are ${n \choose t}$ ways to choose the set
$U$, and at most ${n \choose k}$ ways to choose $k$ vertices outside
$U$.  For each of these vertices, the number of neighbors in $U$ is
distributed as $\bin{t, p}$, and the probability that this Binomial
random variable exceeds $\gamma np$ is at most ${t \choose \gamma np}
p^{\gamma np}$.  Putting this all together, the probability that our
property fails for a certain fixed value of $t$ is at most:
\begin{eqnarray*}
  {n \choose t} \cdot {n \choose k} \cdot \left[
    {t \choose \gamma np} p^{\gamma np}
  \right]^k
  &\leq&
  {n \choose t}^2 \cdot \left[
    {t \choose \gamma np} p^{\gamma np}
  \right]^k \\
  &\leq& \left( \frac{en}{t} \right)^{2t} \cdot \left[
    \frac{etp}{\gamma np}
  \right]^{\gamma np \cdot k} \\
  &=&
  \left( \frac{en}{t} \right)^{2t} \cdot \left[
    \frac{et}{\gamma n}
  \right]^{3t} \\
  &=& \left[ \frac{e^5}{\gamma^3} \cdot \frac{t}{n} \right]^t.
\end{eqnarray*}
Since $t$ is at most $T = \frac{\gamma^3}{2e^5} n$, the final bound is
at most $2^{-t}$.  Summing over all $t$ from $\gamma np$ to $T$, we
see that since $np \rightarrow \infty$, the total failure probability
is still $o(1)$, as desired. \hfill $\Box$

\vspace{3mm}

Now we prove Theorem \ref{thm:fast-lower}(i).  The proof relies on
the following standard definition.

\begin{definition}
  The \textbf{$\boldsymbol{k}$-core} of a graph $G$ is the largest
  induced subgraph with all degrees at least $k$.
\end{definition}

\noindent \textbf{Remark.}\, It is well known that the $k$-core can
always be obtained by repeatedly deleting all vertices of degree less
than $k$, and the result is independent of the order in which these
deletions are performed.

\vspace{3mm}

It is more convenient to prove Theorem \ref{thm:fast-lower}(i) in the
following (equivalent) reparameterized form.

\begin{proposition}
  \label{prop:fast-robber-inverse}
  Let $c > 0$ be fixed, and let $p = \frac{c^3}{30000} n^{c-1}$.  Then
  \whp\ $\Gnp$ has the property that a robber with speed $\frac{1}{c}
  + 2$ can always escape from $n^{1-c}$ cops.
\end{proposition}

\noindent \textbf{Proof.}\, Condition on the high-probability
properties in Lemmas \ref{lem:fast-robber-well-connected},
\ref{lem:fast-robbers-few-edges}, and
\ref{lem:fast-robber-degrees-spread} with $\gamma = c/4$.  Also
condition on the high-probability event that all degrees of $\Gnp$ are
between $0.9 np$ and $1.1 np$.  Note that $np = \frac{c^3}{30000}
n^c$.

Let us specify the robber's winning strategy.  The cops place
themselves first.  Let $C$ be the set of vertices occupied by cops,
and let $C^+$ be the union of $C$ with all immediate neighbors of
vertices in $C$.  Since $|C| \leq n^{1-c}$ and all vertices have
degree at most $1.1 \cdot \frac{c^3}{30000} n^c$, it follows that
$|C^+| \leq \frac{\gamma^3}{2e^5} n$, where $\gamma = c/4$.

Recall that the \emph{$k$-core} of a graph is obtained by repeatedly
deleting all vertices of degree less than $k$.  Let $H$ be the
$\frac{np}{3}$-core of the graph induced by vertices outside $C^+$.
Our first claim is that $H$ always has size at least $(1-c^3) n$.

Indeed, since we conditioned on all degrees exceeding $0.9 np$, as
well as on the result of Lemma \ref{lem:fast-robber-degrees-spread}
with $\gamma = c/4$, the deletion of $C^+$ cannot hurt our minimum
degree condition by very much.  To be precise, the resulting graph has
minimum degree at least $(0.9 - \gamma) np$, except for some small set
of vertices $U_1$ of size at most $\frac{3}{\gamma np} \cdot |C^+| \le
\frac{10^6}{c^4} n^{-c} \cdot |C^+|$.  Applying the same result again,
we find that after deleting $U_1$, the resulting graph has minimum
degree at least $(0.9 - 2\gamma) np$, except for some even smaller set
$U_2$ of size at most $\big( \frac{10^6}{c^4} n^{-c} \big)^2 \cdot
|C^+|$.  Repeatedly applying this result, we see that since $|C^+|
\leq n$, this procedure must certainly terminate within $2/c$
iterations, giving a subgraph with all degrees at least $\big( 0.9 -
\frac{2}{c} \gamma \big) np > \frac{np}{3}$.  The total number of
deleted vertices is at most
\begin{displaymath}
  |C^+| \cdot \left[ 1 + \left( \frac{10^6}{c^4} n^{-c} \right) + \left( \frac{10^6}{c^4} n^{-c} \right)^2 + \cdots \right]
  <
  |C^+| \cdot 2
  < c^3 n,
\end{displaymath}
as claimed.

The robber's strategy is to choose an arbitrary vertex in $H$ for his
initial position.  The cops then make their move, and occupy a new set
of vertices $C'$.  Let $H'$ be the new $\frac{np}{3}$-core of the
graph induced by all vertices except those in $C'$ and its immediate
neighborhood.  Our final claim is that the robber can always move to a
vertex in $H'$.  Clearly, this will imply that the robber can evade
the cops indefinitely.

We must show that there exists a path of length at most $\frac{1}{c} +
2$ from the robber's current position to a vertex in $H'$, which
completely avoids $C'$.  The main observation is that $C' \subset
C^+$, because $C^+$ was defined to include all possible positions of
cops in their next turn.  Therefore, since the robber is in $H$ (the
$\frac{np}{3}$-core of $G \setminus C^+$), he has quite a lot of freedom to
move without running into any of the cops positioned at $C'$.

More precisely, we will show that by traversing at most $\frac{1}{c}
+ 1$ edges in $H$, the robber can already reach $s_0 = \frac{3}{p}
\log n$ vertices.  Indeed, let $S_0=\{x\}, S_1, S_2, \ldots,$ be the
sequence of sets in a breadth-first search performed in $H$ from the
robber's current position $x$.  Let $T_i = S_0 \cup S_1 \cup \cdots
\cup S_i$. Thus $T_{i+1} = T_i\cup N_H(T_i)$.  Now suppose that
$|T_{i+1}|\leq s_0$. It follows from our conditioning on Lemma
\ref{lem:fast-robbers-few-edges} that
\begin{displaymath}
  e(T_{i+1}) \geq |T_i| \left( \frac{np}{3} - 6 \log n \right).
\end{displaymath}
Applying this same result once again we see that
\begin{displaymath}
  |T_{i+1}|
  \geq \frac{e(T_{i+1})}{6\log n}
  \geq \frac{np |T_i|}{20\log n}.
\end{displaymath}
Since $np = \frac{c^3}{30000} n^c$, it follows that if
$i_0=\frac{1}{c}+1$ then $|T_{i_0}| \geq s_0$.

Yet we also conditioned on Lemma \ref{lem:fast-robber-well-connected},
so since $|H'| \geq (1-c^3)n$, there is an edge between $T_{i_0}$ and
$H'$.  Therefore, since the robber is permitted to traverse
$\frac{1}{c} + 2$ edges in a single move, he can indeed land on a
vertex in $H'$ without passing through any vertex (in $C'$) currently
occupied by a cop.  \hfill $\Box$

\section{Concluding remarks}

We have considered the directed version of the classical Cops and
Robbers game, and also the version where the robber moves $R$ edges at
a time, but the cops move only one edge at a time.  Our approach
generalized the best known upper bound to the fast robber setting, but
coincidentally only reproved the same asymptotic in the original
setting.  However, for directed graphs, our general upper bound is
weaker than the corresponding bound for the undirected case.  It would
be nice to obtain an upper bound for directed graphs with asymptotics
similar to our other upper bounds in this paper.

On the other hand, for lower bounds, the fast robber lower bound of
$n^{1-\frac{1}{R-2}}$ we derived is only interesting for $R\geq 5$. It
would be nice to know whether or not an $\omega(\sqrt{n})$ lower bound
can already be achieved for $R=2$.  Another possible version to
address is when the cops and the robber both move at the \emph{same}\/
speed $R > 1$. Our upper bound on the number of cops in the fast
robber scenario still carries over, since faster cops are more
powerful. It would be interesting to decide whether there is a better
lower bound of $\omega(\sqrt{n})$ for this case.


\appendix

\section{Routine inequalities}

\begin{inequality}
  \label{ineq:di-p}
  Let $p = \frac{13 (\log \log n)^2}{\log n}$, and let $r = \frac{6}{p}
  \log \frac{4}{p}$.  Then
  \begin{displaymath}
   2 p^{-r} (1+p)^r \leq \frac{np}{2}
    \andsp
    \text{and}
    \andsp
    (1+p)^r \geq \frac{16}{p} \log n.
  \end{displaymath}
\end{inequality}

\noindent \textbf{Proof.}\, The left hand side of the first
inequality is at most
\begin{displaymath}
 2 p^{-r} (1+p)^r
  \leq
 2 p^{-r} e^{pr}
  = 2 p^{-r} \left( \frac{4}{p} \right)^6
  \leq
  \left( \frac{4}{p} \right)^{r+6}.
\end{displaymath}
Therefore, it suffices to show that $\big(\frac{4}{p}\big)^{r+7} \leq
n$, or equivalently, that
\begin{displaymath}
  (r+7) \log \frac{4}{p} \leq \log n
\end{displaymath}
Yet $r+7 \leq 2r$, and $2r \log \frac{4}{p} = \frac{12}{p} \big( \log
\frac{4}{p} \big)^2$, so plugging in the definition of $p$, we see
that this is indeed less than $\log n$.

For the second inequality, since $p$ is small,
\begin{displaymath}
  (1+p)^r
  \geq
  e^{pr/2}
  =
  \left( \frac{4}{p} \right)^3
  >
  \frac{16}{p} \log n,
\end{displaymath}
since $\frac{1}{p} = \frac{\log n}{13 (\log \log n)^2}$.  This
completes the proof.  \hfill $\Box$

\vspace{3mm}

\begin{inequality}
  \label{ineq:p}
  There is a function $p = p(n) = 2^{-(1-o(1)) \sqrt{\log_2 n}}$ and a
  positive integer $l = (1+o(1)) \sqrt{\log_2 n}$ such that when we
  define $k = \frac{16}{p} \log n$, we have the inequalities
  \begin{displaymath}
    k^{l+1} (1+p)^{2^l} \leq np/2
    \andsp
    \text{and}
    \andsp
    (1+p)^{2^l} \geq k.
  \end{displaymath}
\end{inequality}

\noindent \textbf{Proof.}\, First observe that we will have $\log
\frac{1}{p} = \Theta(\sqrt{\log n})$, so $\log k = (1+o(1)) \log
\frac{1}{p}$.  Let $l$ be the smallest positive integer for which the
second inequality is satisfied.  This immediately gives $(1+p)^{2^l}
\leq k^2$.  Also,
\begin{displaymath}
  l = \left\lceil
  \log_2 \left( \frac{\log k}{\log (1+p)} \right)
  \right\rceil
  =
  \left\lceil
    \log_2 \left( (1+o(1)) \frac{1}{p} \log \frac{1}{p} \right)
  \right\rceil
  =
  (1+o(1)) \log_2 \frac{1}{p}.
\end{displaymath}
To establish the first inequality, we have $k^{l+1} (1+p)^{2^l} \leq
k^{l+1} k^2$, so it suffices to show that
\begin{equation}
  \label{ineq:p-internal}
  (l+3) \log_2 k \leq \log_2 \frac{np}{2}
\end{equation}
From the asymptotics of $p$, we have
\begin{eqnarray*}
  l+3 &=& (1+o(1)) \log_2 \frac{1}{p} \\
  \log_2 k &=& (1+o(1)) \log_2 \frac{1}{p} \\
  \log_2 \frac{np}{2} &=& (1-o(1)) \log_2 n,
\end{eqnarray*}
so it is clear that Inequality \eqref{ineq:p-internal} is satisfied
for an appropriate choice of $p = 2^{-(1-o(1)) \sqrt{\log_2 n}}$.
\hfill $\Box$

\vspace{3mm}

\begin{inequality}
  \label{ineq:fast-p}
  Let $R > 1$ be given, and define the sequence $d_0, d_1, \ldots$ via
  the recursion $d_0 = R$, $d_{i+1} = d_i + \big\lceil \frac{d_i}{R}
  \big\rceil$.  Let $r_i = \big\lceil \frac{d_i}{R} \big\rceil$.  Then
  there is a function
  \begin{displaymath}
    p =
    p(n) =
    \left( 1 + \frac{1}{R} \right)^{-(1-o(1)) \sqrt{\log_{1 + \frac{1}{R}} n}}
  \end{displaymath}
  and a positive integer $l = (1+o(1)) \sqrt{\log_{1 + \frac{1}{R}}
    n}$ such that when we define $k = \frac{16}{p} \log n$, we have
  the inequalities
  \begin{displaymath}
    k^{2^R} \cdot k^l (1+p)^{r_l} \leq \frac{np}{2}
    \andsp
    \text{and}
    \andsp
    (1+p)^{r_l} \geq k.
  \end{displaymath}
\end{inequality}

\noindent \textbf{Proof.}\, The proof is nearly identical to the
previous lemma.  We will have $\log \frac{1}{p} = \Theta(\sqrt{\log
  n})$, so $\log k = (1+o(1)) \log \frac{1}{p}$.  Let $l$ be the
smallest positive integer for which the second inequality is
satisfied.  Since $r_{l+1} \leq 2 r_l$, this immediately gives
$(1+p)^{r_l} \leq k^2$.

Let us estimate an asymptotic upper bound for $l$.  Observe that $d_i
\geq \big( 1 + \frac{1}{R} \big) d_{i-1}$, so $r_l \geq \big( 1 +
\frac{1}{R} \big)^l$.  Hence if we let $l'$ satisfy
\begin{displaymath}
  (1+p)^{ (1 + \frac{1}{R} )^{l'} } = k,
\end{displaymath}
then $l \leq l'$.  Yet
\begin{displaymath}
  l' =
  \log_{1 + \frac{1}{R}} \frac{\log k}{\log (1+p)}
  =
  \log_{1 + \frac{1}{R}} \left( (1+o(1)) \frac{1}{p} \log \frac{1}{p} \right)
  =
  (1+o(1)) \log_{1 + \frac{1}{R}} \frac{1}{p}.
\end{displaymath}

To establish the first inequality, we initially noted that
$(1+p)^{r_l} \leq k^2$, so $k^{2^R} \cdot k^l (1+p)^{r_l} \leq k^{2^R
  + 2 + l}$.  Thus it suffices to show that
\begin{equation}
  \label{ineq:fast-p-internal}
  (2^R + 2 + l) \log_{1 + \frac{1}{R}} k \leq \log_{1 + \frac{1}{R}} \frac{np}{2}
\end{equation}
From the asymptotics of $p$, we have
\begin{eqnarray*}
  2^R + 2 + l &\leq& (1+o(1)) l'
  \ \ = \ \ (1+o(1)) \log_{1 + \frac{1}{R}} \frac{1}{p} \\
  \log_{1 + \frac{1}{R}} k &=& (1+o(1)) \log_{1 + \frac{1}{R}} \frac{1}{p} \\
  \log_{1 + \frac{1}{R}} \frac{np}{2} &=& (1-o(1)) \log_{1 + \frac{1}{R}} n,
\end{eqnarray*}
so inequality \eqref{ineq:fast-p-internal} is clearly satisfied by
appropriately choosing $\log_{1 + \frac{1}{R}} \frac{1}{p} = (1-o(1))
\sqrt{\log_{1 + \frac{1}{R}} n}$.  This is precisely the asymptotic
claimed in our statement.  \hfill $\Box$

\end{document}